\def \NN {\mathbb N}
\def \CC {\mathbb C}

\def \RR {\mathbb R}
\def \ZZ {\mathbb Z}

\def \epsilon{\varepsilon}

\def \K  {{\mathcal K}}

\def \LL  {{\mathcal L}}
\def \P  {{\mathcal P}}

\def \R  {{\mathcal R}}

\def \d {\text{d}}

\def \fine {{\hfill \qedsymbol}}

\def \Ga {\Gamma}
\def \si {\sigma}

\newcommand{\res}{\text{res}}
\renewcommand{\S}{{\mathcal S}}

\documentclass[12pt,reqno]{amsart}

\usepackage{amsmath,amssymb}

\usepackage{url} 

\numberwithin{equation}{section}

\begin{document}


\title[]{On the standard twist of the $L$-functions \\ of half-integral weight cusp forms}
\author[]{J.KACZOROWSKI \lowercase{and} A.PERELLI}
\maketitle

{\bf Abstract.} The standard twist $F(s,\alpha)$ of $L$-functions $F(s)$ in the Selberg class has several interesting properties and plays a central role in the Selberg class theory. It is therefore natural to study its finer analytic properties, for example the functional equation. Here we deal with a special case, where $F(s)$ satisfies a functional equation with the same $\Ga$-factor of the $L$-functions associated with the cusp forms of half-integral weight; for simplicity we present our results directly for such $L$-functions. We show that the standard twist $F(s,\alpha)$ satisfies a functional equation reflecting $s$ to $1-s$, whose shape is not far from a Riemann-type functional equation of degree 2 and may be regarded as a degree 2 analog of the Hurwitz-Lerch functional equation. We also deduce some result on the growth on vertical strips and on the distribution of zeros of $F(s,\alpha)$.

\smallskip
{\bf Mathematics Subject Classification (2000):} 11M41, 11F66, 11F37

\smallskip
{\bf Keywords:} Selberg class, standard twist, functional equations, half-integral weight modular forms

\vskip.5cm
\section{Introduction}

\smallskip
{\bf 1.1.~Motivations.} In \cite{Ka-Pe/1999a}, \cite{Ka-Pe/2005} and \cite{Ka-Pe/resoI} we introduced and studied the standard twist
\begin{equation}
\label{1-1}
F(s,\alpha) = \sum_{n=1}^\infty \frac{a(n)}{n^s} e(-\alpha n^{1/d})
\end{equation}
of any given function $F(s)$ from the extended Selberg class $\S^\sharp$; here $a(n)$ and $d>0$ are, respectively, the Dirichlet coefficients and the degree of $F(s)$, $e(x)=e^{2\pi ix}$ and $\alpha>0$. We refer to the next subsection for precise definitions of the quantities we introduce in this subsection; here we recall only that the class $\S^\sharp$ consists, roughly, of the Dirichlet series with analytic continuation to the whole complex plane and satisfying a general functional equation of Riemann type. In particular, essentially all $L$-functions arising from number theory and automorphic forms belong, at least conditionally, to $\S^\sharp$.

\smallskip
The main known properties of $F(s,\alpha)$ are as follows; in order to keep a lighter notation, we shall always assume that $F(s)$ is normalized in a way to be explained in Subsection 1.2; such a normalization is always satisfied by the classical $L$-functions. For every $F\in\S^\sharp$ and $\alpha>0$, the series in \eqref{1-1} converges absolutely for $\si>1$ and $F(s,\alpha)$ has meromorphic continuation to $\CC$. More precisely, if $\alpha$ does not belong to the spectrum Spec$(F)$ of $F(s)$, an infinite discrete subset of $\RR^+$ to be defined in the next subsection, then $F(s,\alpha)$ is an entire function of finite order, while if $\alpha \in$ Spec$(F)$, then its only singularities are at most simple poles at the points
\begin{equation}
\label{1-1-1}
s_\ell= \frac{d+1}{2d} -\frac{\ell}{d}, \quad \ell=0,1,\dots
\end{equation}
In addition,
\[
\res_{s=s_0} F(s,\alpha) = C_F \frac{\overline{a(n_\alpha)}}{n_\alpha^{1-s_0}}, \qquad C_F\neq 0;
\]
hence $F(s,\alpha)$ has always a simple pole at $s=s_0$ when $\alpha \in$ Spec$(F)$. Moreover, $F(s,\alpha)$ has 
\eject
\noindent
polynomial growth on every vertical strip, although the known bounds are weak in general. We refer to our papers \cite{Ka-Pe/2005} and \cite{Ka-Pe/resoI} for these and other results.

\smallskip
At present, the interest of the standard twist comes mainly from the fact that it plays a central role in the Selberg class theory. The main aim of such a theory is describing the structure of the Selberg class $\S$, roughly the subclass of the functions $F\in\S^\sharp$ with a general Euler product and satisfying the Ramanujan conjecture $a(n) \ll n^\epsilon$. It is expected that $\S$ coincides with the class of automorphic $L$-functions and hence, in particular, that the degree $d$ is always an integer. The above properties of the standard twist, often coupled with those of other nonlinear twists, have been exploited to verify such an expectation for every degree in the range $0<d<2$; see Conrey-Ghosh \cite{Co-Gh/1993} and our papers \cite{Ka-Pe/1999a} and \cite{Ka-Pe/2011}. A very simple example is the following proof of the nonexistence of functions in $\S^\sharp$ with degree $0<d<1$: for $d$ in such a range we have that $s_0>1$, see \eqref{1-1-1}, hence $F(s,\alpha)$ would have a pole in the region of absolute convergence.

\smallskip
It is therefore natural to ask for a description of the finer analytic properties of $F(s,\alpha)$, and in particular to rise the following problems.

\smallskip
(i) Does $F(s,\alpha)$ satisfy a functional equation relating $s$ to $1-s$ ? We refer to the discussion at the end of Subsection 1.3 for additional information related to this problem.

(ii) Study of the finer polar structure of $F(s,\alpha)$, in particular the existence of finitely or infinitely many poles at the points \eqref{1-1-1}; examples of both type exist, see Remark 3 below.

(iii) Give precise convexity bounds for the Lindel\"of $\mu$-function of $F(s,\alpha)$
\begin{equation}
\label{1-1-3}
\mu(\si) =\mu_F(\si,\alpha) =  \inf\{\xi: F(\si+it,\alpha) \ll (1+|t|)^\xi \ \text{as} \ |t|\to\infty\}.
\end{equation}

(iv) Determine location and counting of the zeros of $F(s,\alpha)$, distinguishing between trivial and nontrivial zeros.

(v) Other analytic problems on $F(s,\alpha)$ like bounds for moments, sharp uniform bounds in $\alpha$, etc.

\smallskip
\noindent
None of these problems is solved at present when the degree $d$ of $F(s)$ is $\geq2$; in this paper we give a first contribution, in the framework of a special family of degree 2 $L$-functions.

\medskip
{\bf 1.2.~Definitions and notation.} The extended Selberg class $\S^\sharp$ consists of non identically vanishing Dirichlet series $F(s)$, absolutely convergent for $\si>1$, such that $(s-1)^mF(s)$ is entire of finite order for some integer $m$ and satisfying a functional equation of type
\begin{equation}
\label{1-2}
Q^s\prod_{j=1}^r\Gamma(\lambda_js+\mu_j) F(s) = \omega Q^{1-s}\prod_{j=1}^r\Gamma(\lambda_j(1-s)+\overline{\mu_j}) \overline{F(1-\overline{s})}
\end{equation}
with $|\omega|=1$, $Q>0$, $\lambda_j>0$ and $\Re{\mu_j}\geq0$. We refer to Selberg \cite{Sel/1989}, Conrey-Ghosh \cite{Co-Gh/1993}, to our survey papers \cite{Ka-Pe/1999b}, \cite{Kac/2006}, \cite{Per/2005}, \cite{Per/2004}, \cite{Per/2010}, \cite{Per/2017} and to the forthcoming book \cite{Ka-Pe/book} for definitions, examples and the basic theory of the class $\S^\sharp$. We recall that degree $d$, conductor $q$ and $\xi$-invariant of $F(s)$ are defined, respectively, by
\begin{equation}
\label{1-3}
d=2\sum_{j=1}^r\lambda_j, \quad q= (2\pi)^dQ^2\prod_{j=1}^r\lambda_j^{2\lambda_j}, \quad \xi = 2\sum_{j=1}^r (\mu_j-\frac12) = \eta+id \theta
\end{equation}
with $\eta,\theta\in\RR$. Throughout the paper we always assume that $F(s)$ is normalized by the condition $\theta=0$.
Moreover, we write $n_\alpha=q d^{-d}\alpha^d$ and the spectrum of $F(s)$ is defined as
\[
\text{Spec$(F) = \{\alpha>0: n_\alpha\in\NN \ \text{with} \ a(n_\alpha)\neq 0\}$}.
\]

\medskip
{\bf 1.3.~Set up of the problem.}  As already pointed out, to gain experience on the problems listed in Subsection 1.1 we consider the following special case related to half-integral weight modular forms. Their Hecke $L$-functions satisfy a special functional equation, and by a lucky coincidence we can apply certain methods developed after Linnik in our papers, and in particular in  \cite{Ka-Pe/2017}.

Let $f$ be a cusp form of half-integral weight $\kappa=k/2$ and level $N$, where $k>0$ is an odd integer and $4|N$, and $L_f(s)$ be the associated Hecke $L$-function. Then $L_f(s)$ is entire and satisfies the functional equation
\begin{equation}
\label{1-4}
\Lambda_f(s) = \omega\Lambda_{f^*}(\kappa-s), \qquad \text{where} \qquad \Lambda_f(s) = \big(\frac{\sqrt{N}}{2\pi}\big)^s\Gamma(s)L_f(s),
\end{equation}
$\omega=i^{-\kappa}$ and $f^*$ is related to $f$ by $f^*(z) = (\sqrt{N}z)^{-\kappa} f(-1/Nz)$. Note that $L_{f^*}(s)$ is also entire and has properties similar to $L_f(s)$. We refer to Ogg \cite{Ogg/1969} and to Section 4.3 of Miyake \cite{Miy/1989} for the basic analytic theory of modular forms, and to Bruinier \cite{Bru/1997} for a detailed exposition of the half-integral weight case.

\smallskip
Comparing functional equations \eqref{1-2} and \eqref{1-4}, it is clear that the function $L_f(s)$ does not belong to the extended Selberg class $\S^\sharp$. However, it comes close after the normalization $s\mapsto s+\frac{\kappa-1}{2}$. Indeed, writing
\begin{equation}
\label{1-5}
F(s) = L_f(s+\frac{\kappa-1}{2}) \quad \text{and} \quad F^*(s) = L_{f^*}(s+\frac{\kappa-1}{2}),
\end{equation}
the functional equation \eqref{1-4} becomes
\begin{equation}
\label{1-6}
\big(\frac{\sqrt{N}}{2\pi}\big)^s\Gamma(s + \frac{\kappa-1}{2})F(s) = \omega\big(\frac{\sqrt{N}}{2\pi}\big)^{1-s}\Gamma(1-s + \frac{\kappa-1}{2})F^*(1-s).
\end{equation}
Although this is not exactly of the form \eqref{1-2} (indeed $F^*(s)$ is not necessarily equal to $\overline{F(\overline{s})}$, and for $k=1$ we have $\frac{\kappa-1}{2}=\frac{k-2}{4}<0$), most results in the Selberg class theory hold in this case as well. Note that the other requirements of $\S^\sharp$ are satisfied by $F(s)$, as it can be checked by the argument in the proof of Theorem 5.1 and Corollary 5.2 of Iwaniec \cite{Iwa/1997} (see also the proof of the Theorem in Kaczorowski {\it et al.} \cite{KMPSW/2006} and on p.217-218 of Carletti {\it et al.} \cite{C-M-P/2009}). In particular, both $F(s)$ and $F^*(s)$ are absolutely convergent for $\sigma>1$. Note also that, in view of \eqref{1-6}, $F(s)$ has degree $d=2$ and conductor $q=N$.

\smallskip
From now on we consider cusp forms $f$ of weight $\kappa$ and level $N$ as above, and the normalized Hecke $L$-functions $F(s)$ and $F^*(s)$ as in \eqref{1-5}. We denote by $a(n)$ and $a^*(n)$ their Dirichlet coefficients, respectively, and let $F(s,\alpha)$ be the standard twist of $F(s)$, defined by \eqref{1-1} with $d=2$ in this case. Then, minor formal modifications to the proof of Theorem 1 of \cite{Ka-Pe/2005}, or of Theorems 1 and 2 of \cite{Ka-Pe/resoI}, give the following properties of $F(s,\alpha)$. Writing in this special case
\begin{equation}
\label{1-7}
n_\alpha = N\alpha^2/4 \qquad \text{and} \qquad \text{Spec$^*(F) = \{\alpha>0: a^*(n_\alpha)\neq0\}$},
\end{equation}
where $a^*(n_\alpha)=0$ if $n_\alpha\not\in\NN$, we have that $F(s,\alpha)$ is entire if $\alpha\not\in$ Spec$^*(F)$. If $\alpha\in$ Spec$^*(F)$, then $F(s,\alpha)$ is meromorphic over $\CC$ with at most simple poles at
\begin{equation}
\label{1-8}
s_\ell = \frac{3}{4} - \frac{\ell}{2} \hskip1.5cm \ell = 0,1,\dots
\end{equation}
(in this case the condition $\theta=0$ is satisfied) and
\begin{equation}
\label{1-9}
\res_{s=s_0} F(s,\alpha) = C_F \frac{a^*(n_\alpha)}{n_\alpha^{1/4}} \hskip1.5cm C_F\neq 0.
\end{equation}
Moreover, in all cases $F(s,\alpha)$ has polynomial growth on vertical strips. 

\smallskip
The main result of this paper is the affirmative answer to problem (i) in Subsection 1.1 in this special case, i.e. $F(s,\alpha)$ satisfies a functional equation relating $s$ to $1-s$ by means of the $\Ga$ function. Once this is obtained, more or less standard techniques allow to solve some of the other problems in the list, in particular (ii), (iii) and (iv). However, the output is somehow unconventional in some cases, e.g. for the trivial zeros. Of course, the functional equation of $F(s,\alpha)$ may be of independent interest inside the modular forms theory, as it provides new information about the associated $L$-functions.

\smallskip
Before stating the main theorem, we recall that in \cite{Ka-Pe/2014a} we proved that if a function $F\in\S^\sharp$ has degree $\geq2$ and satisfies the Ramanujan conjecture (i.e. $a(n)\ll n^\epsilon$), then the standard twist $F(s,\alpha)$ does not belong to $\S^\sharp$. Actually, a variant of the arguments in \cite{Ka-Pe/2014a}, which we sketch in the Appendix below, proves under the same assumptions the stronger assertion that $F(s,\alpha)$ does not satisfy a functional equation of type \eqref{1-2}. Nevertheless, in the case of half-integral weight modular forms, $F(s,\alpha)$ satisfies a functional equation reflecting $s$ into $1-s$ via suitable $\Ga$-factors.

\medskip
{\bf Remark 1.} The Ramanujan conjecture is crucial for our results in \cite{Ka-Pe/2014a}. Indeed, for certain Dirichlet $L$-functions, say $L_1(s)$, the standard twists $L_1(s,\alpha)$ are still degree 1 functions of $\S^\sharp$ for suitable values of $\alpha$. But such $L_1(s)$ can be lifted to degree 2 functions $L_2(s)$ by letting $s\mapsto 2s-1/2$, and for the same values of $\alpha$ their standard twists $L_2(s,\alpha)$ satisfy a degree 2 functional equation of type \eqref{1-6}, with $\kappa =1/2$ or $\kappa = 3/2$. Moreover, when $\kappa=3/2$ both $L_2(s)$ and $L_2(s,\alpha)$ belong to $\S^\sharp$. However, such $L_2(s)$'s do not satisfy the Ramanujan conjecture, thus showing that this hypothesis is crucial in \cite{Ka-Pe/2014a}; also, it turns out that such $L$-functions are related with half-integral weight modular forms. We refer to Section 2 of \cite{Ka-Pe/abs} for the basic properties of lifts in $\S^\sharp$. Apart from these special cases, we don't know examples of functions $F(s)$, satisfying a functional equation of type \eqref{1-2} with degree $d\geq 2$ but not the Ramanujan conjecture, with standard twist $F(s,\alpha)$ also satisfying a functional equation of type \eqref{1-2}. \fine

\medskip
{\bf 1.4.~Functional equation.} Let $f$ be a cusp form of half-integral weight $\kappa=k/2$ and level $N$, $F(s)$ and $F^*(s)$ be as in  \eqref{1-5}, $F(s,\alpha)$ be the standard twist of $F(s)$, and $a^*(n)$ be the coefficients of $F^*(s)$. We write
\begin{equation}
\label{1-10}
k=2h+1 \ \text{with}\  h\in\{0,1,2\dots\}, \ \ \mu= (2h-1)/4, \ \ h^*=\max(0,h-1), 
\end{equation}
\begin{equation}
\label{1-11}
\text{$\nu=\pm\sqrt{n}$ \ with $n=1,2,\dots$ and \ $\nu_\alpha=\sqrt{n_\alpha}$ \ with $n_\alpha$ as in \eqref{1-7}.}
\end{equation}
Moreover, for $\ell=0,\dots,h^*$ we also write
\begin{equation}
\label{1-12}
F^*_\ell(s,\alpha) = e^{-i\pi s} F^+_\ell(s,\alpha) + e^{i\pi s} F^-_\ell(s,\alpha),
\end{equation}
where, putting
\begin{equation}
\label{1-13}
c^*(\nu^2) = c_\ell^*(\nu^2) =
\begin{cases}
-e^{i\pi\mu} a^*(\nu^2) \ \ &\text{if} \ \  \nu\geq 1 \\
e^{i\pi(1/2+\ell-\mu)} a^*(\nu^2) &\text{if} \ \ -\nu_\alpha < \nu\leq -1 \\
e^{-i\pi\mu} a^*(\nu^2) \ \ &\text{if} \ \  \nu<-\nu_\alpha,
\end{cases}
\end{equation}
the generalized Dirichlet series $F^\pm_\ell(s,\alpha)$ are defined by
\begin{equation}
\label{1-14}
F^+_\ell(s,\alpha) = \sum_{\nu>-\nu_\alpha} \frac{c^*(\nu^2)}{|\nu|^{1/2+\ell} |\nu+\nu_\alpha|^{2s-1/2-\ell}}, \quad F^-_\ell(s,\alpha) = \sum_{\nu<-\nu_\alpha} \frac{c^*(\nu^2)}{|\nu|^{1/2+\ell} |\nu+\nu_\alpha|^{2s-1/2-\ell}};
\end{equation}
see also \eqref{2.25} for more information on the shape of \eqref{1-14}. Note that, thanks to the convergence properties of $F^*(s)$, such generalized Dirichlet series are absolutely convergent for $\si>1$. Note also that the second range in \eqref{1-13} is empty if $0\leq n_\alpha\leq 1$. Finally we define the coefficients $a_\ell=a_\ell(h^*)$ by means of the following polynomial identity
\begin{equation}
\label{1-15}
\prod_{1\leq j \leq h^*}(X+2j-1) = \sum_{\ell=0}^{h^*} a_\ell \prod_{0\leq \nu \leq h^*-1-\ell} (X+\nu),
\end{equation}
where, throughout the paper, an empty product equals 1; note that $a_0=1$ for every $h^*\geq0$.  With the above notation we have

\medskip
{\bf Theorem.} {\sl Let $\alpha>0$ and $\ell=0,\dots,h^*$. Then the functions $F^*_\ell(s,\alpha)$ are entire and $F(s,\alpha)$ satisfies the functional equation}
\begin{equation}
\label{1-16}
F(s,\alpha) = \frac{\omega}{i\sqrt{2\pi}} \big(\frac{\sqrt{N}}{4\pi}\big)^{1-2s} \sum_{\ell=0}^{h^*} a_\ell \Gamma\big(2(1-s)-1/2- \ell\big) F^*_\ell(1-s,\alpha).
\end{equation}

\medskip
Note that functional equation \eqref{1-16} is not exactly of Riemann type, but not far from it. Indeed, one can see from the proofs of Lemma 2.1 and Corollary 1 below that the functions $F^\pm_\ell(s,\alpha)$ can be expressed as a kind of stratification of $F^*(s)$; by this we mean that each $F^\pm_\ell(s,\alpha)$ is related to a combination of shifts of $F^*(s)$, see \eqref{2.29}, \eqref{2-37} and \eqref{2-38}. Therefore, \eqref{1-16} resembles the asymmetric form of the functional equation of $F(s)$, and hence $F(s,\alpha)$ is expected to behave essentially as a degree 2 function in the extended Selberg class $\S^\sharp$; this is confirmed by the results below. Note also that letting $\alpha=0$ in \eqref{1-16} one obtains the asymmetric form of the functional equation of $F(s)$. Indeed, for every $\ell=0,\dots,h^*$ we have
\[
F^*_\ell(1-s,0) = \big(e^{i\pi(s+\mu)}- e^{-i\pi(s+\mu)}\big)F^*(1-s),
\]
and following the initial steps of the proof of the Theorem one can easily rebuild the $\Gamma$-factor of $F(s)$ from the $\Gamma$-factors in \eqref{1-16}. 

\medskip
{\bf Remark 2.} We obtain the functional equation in the theorem thanks to the special form of the $\Gamma$-factor in \eqref{1-6}, which enables the explicit computation of a certain hypergeometric function arising in the proof, see \eqref{2.11} and \eqref{2.12} below. Obviously, we can get the same result for any $F\in\S^\sharp$ with the same $\Ga$-factor in its functional equation. Moreover, similar arguments can be carried over for a certain class of $\Gamma$-factors, for example of type $\Gamma(ds/2 + \mu)$ with any $d\geq 1$ and suitable real values of $\mu$, thus producing analogous functional equations for the standard twists of the corresponding $L$-functions. Here we restrict ourselves to the case of \eqref{1-5} for simplicity, since it already covers an interesting classical case. We also remark that, as it is, our method does not apply to the case of classical modular forms of integral weight, nor to the Maass forms. Actually, we don't even expect a functional equation as in the theorem in such cases. For example, we know from \cite{Kac/2004} that there are integral weight modular forms whose associated standard twist has infinitely many poles at the points $s_\ell$ in \eqref{1-1-1} when $\alpha$ is in the spectrum, and this is not coherent with a functional equation of type \eqref{1-16}. \fine

\medskip
{\bf Remark 3.} Let $F\in\S^\sharp$ and $\alpha\in$ Spec$(F)$. We already pointed out in Subsection 1.1 that $F(s,\alpha)$ has at most simple poles at the points $s=s_\ell$ in \eqref{1-1-1}, and from Remark 2 we see that there cases where $s=s_\ell$ is actually a pole for infinitely many $\ell$'s. On the other hand, only $s=s_0$ is a pole in the case $F(s)=\zeta(s)$. The existence of infinitely or finitely many poles at the points $s=s_\ell$ depends on the structure of certain quotients of $\Gamma$ functions related to the hypergeometric function arising as Mellin transform of the $\Gamma$-factors in the functional equation of $F(s)$. In the case treated in the theorem, it follows directly from the functional equation \eqref{1-16} that $F(s,\alpha)$ is holomorphic apart possibly at $s=s_\ell$ with $\ell=0,\dots,h^*$, since the functions $F^*_\ell(s,\alpha)$ are entire. An explicit expression for the residue $\kappa_\ell(\alpha)$ of $F(s,\alpha)$ at $s=s_\ell$, $\ell=0,\dots,h^*$, is given in \eqref{2.21} below. Moreover, \eqref{1-16} provides an alternative expression for such residues. We also note that, whatever the value of $\alpha$ is, the coefficient $a^*(n_\alpha)$ never appears on the right hand side of \eqref{1-16} (see \eqref{1-12} and \eqref{1-14}), but it pops up on the left hand side of \eqref{1-16} when $\alpha\in$ Spec$^*(F)$, inside the residue at $s=s_0$ (see \eqref{1-9}).  \fine

\medskip
In the next three subsections we prove some corollaries of our main result; hence $F(s)$ is as in \eqref{1-5} and $F(s,\alpha)$ is its standard twist in such corollaries. We also introduce several constants, sometimes explicitly depending on $\alpha$; it is however clear that in general all such constants may depend on $f$ and $\alpha$.

\medskip
{\bf 1.5.~Order of growth.} Functional equation \eqref{1-16} and the properties of the functions $F^*_\ell(s,\alpha)$ open the possibility of a further study of the standard twist $F(s,\alpha)$. We already know that $F(s,\alpha)$ has polynomial growth on vertical strips, see Theorem 2 of \cite{Ka-Pe/2005}, but the bounds there are definitely weak from a quantitative viewpoint. The functions $F_\ell^\pm(s,\alpha)$ have polynomial growth on vertical strips as well, see Lemma 2.1 below, and for each choice of $\pm$ we denote by
\begin{equation}
\label{1-17}
\mu^\pm(\si)=\inf\{\xi: F_0^\pm(\si+it,\alpha) \ll |t|^\xi \ \text{as} \ t\to\pm\infty\}
\end{equation}
the one-sided Lindel\"of $\mu$-function of $F^\pm_0(s,\alpha)$. One checks by standard arguments (see Section 9.41 of Titchmarsh \cite{Tit/1939}) that the main properties of the Lindel\"of $\mu$-function of Dirichlet series are satisfied by $\mu^\pm(\si)$ as well, namely $\mu^\pm(\si)$ is continuous, convex, non-negative and strictly decreasing until it becomes identically vanishing. Moreover, let
\begin{equation}
\label{1-18}
\mu^*(\si) = \max\big(\mu^+(\si), \mu^-(\si)\big).
\end{equation}
Under the same assumptions of the previous subsection, with the above notation and recalling \eqref{1-1-3} we have the following degree 2 convexity bound for $F(s,\alpha)$.

\medskip
{\bf Corollary 1.} {\sl Let $\alpha>0$ and $\ell=0,\dots,h^*$. Then the functions $F^\pm_\ell(s,\alpha)$ are entire with polynomial growth on vertical strips, and}
\[
\mu(\si)= 1-2\si+\mu^*(1-\si).
\]

\medskip
Hence
\[
\mu(\si)=0 \ \ \text{for} \ \ \si\geq1 \ \ \text{and} \ \ \mu(\si)=1-2\si \ \ \text{for} \ \ \si\leq 0,
\]
and the same holds for $\mu^*(\si)$. Comparing with the bounds for the general standard twist in Theorem 2 of \cite{Ka-Pe/2005}, we see that Corollary 1 is definitely sharper, of course in the special case under consideration.

\medskip
{\bf 1.6.~Trivial zeros.} Our next two results study the zeros of $F(s,\alpha)$, denoted as usual by $\rho=\beta+i\gamma$. We first study the trivial zeros, although distinguishing between trivial and nontrivial zeros is more delicate and unconventional in this case. Indeed, due to the shape of $F_\ell^*(s,\alpha)$ and to the lack of Euler product, trivial zeros come from the interferences between the two terms on the right hand side of \eqref{1-12}, rather than from the poles of the $\Gamma$-factors as in the classical cases. As it will be clear in a moment, in this paper we restrict ourselves to a rough definition of trivial zeros, which in principle does not locate them uniquely. Recalling \eqref{1-10}-\eqref{1-13} let
\[
\begin{split}
&\text{$\nu_+=\nu_+(\alpha)$  be the value $\nu>-\nu_\alpha$  for which $c^*(\nu^2)\neq0$ and $|\nu+\nu_\alpha|$ is minimum,} \\
&\text{$\nu_-=\nu_-(\alpha)$  be the value $\nu<-\nu_\alpha$  for which $c^*(\nu^2)\neq0$ and $|\nu+\nu_\alpha|$ is minimum}
\end{split}
\]
and write
\begin{equation}
\label{1-19}
m_\pm = m_\pm(\alpha) = |\nu_\pm + \nu_\alpha| \quad \text{and} \quad c^*_\pm  = \sqrt{m_\pm} \frac{c_0^*(\nu_\pm^2)}{|\nu_\pm|^{1/2}} = \rho_\pm e^{i\theta_\pm};
\end{equation}
clearly, $m_\pm,\rho_\pm> 0$, $\theta_\pm\in[0,2\pi)$ and $\rho_\pm,\theta_\pm$ depend on $\alpha$. Further, writing as usual $s=\sigma+it$, let $\ell(\alpha)$ be the line of equation
\[
 t=\frac{\si}{\pi}\log\left(\frac{m_+}{m_-}\right) +\frac{1}{2\pi}\log\left(\frac{\rho_+m_-^2}{\rho_-m_+^2}\right),
\]
and for $\epsilon>0$ let
\begin{equation}
\label{1-20}
\LL_\epsilon(\alpha) = \left\{\text{$s\in\CC$ with distance $<\epsilon$ from the line} \ \ell(\alpha)\right\}.
\end{equation}
Then the trivial zeros of $F(s,\alpha)$ may be defined as the zeros contained in the part of the tubular region $\LL_\delta(\alpha)$ with $\si<-\si^-$, with suitable $\delta>0$ and $\si^-\geq 0$. Indeed, under the same assumptions of Subsection 1.4 we have

\medskip
{\bf Corollary 2.} {\sl Let $\alpha>0$. Then there exist $\delta>0$ and $\si^-\geq 0$ such that $F(s,\alpha)\neq 0$ for $\si<-\si^-$ unless $s\in\LL_\delta(\alpha)$. More precisely, there exists $c_1(\alpha)>0$ with the property that for every $0<\epsilon<c_1(\alpha)$  there exists $\si_\epsilon\geq0$ such that 

i) $F(s,\alpha)\neq 0$ for $\si<-\si_\epsilon$ unless $s\in\LL_{\epsilon}(\alpha)$,

ii) there exist infinitely many zeros of $F(s,\alpha)$ in $\LL_{\epsilon}(\alpha)$ with $\beta<-\si_\epsilon$.}

\medskip
Hence we may choose as $\si^-$ any fixed value larger than the infimum of the $\si_\epsilon$'s, and $\delta$ accordingly, although in principle this does not define uniquely the trivial zeros. However, denoting by $T_\epsilon(R, \alpha)$ the number of the zeros in $\LL_\epsilon(\alpha)$ with $-R\leq \beta<-\si_\epsilon$, it is clear from the proof of Corollary 2 that as $R\to+\infty$
\begin{equation}
\label{1-21}
T_\epsilon(R,\alpha) = c_2(\alpha) R +O_\epsilon (1)
\end{equation}
for any given $\alpha>0$, with some $c_2(\alpha)\neq0$ (see in particular \eqref{3-60} below). Hence trivial zeros are well defined by the above choice, apart from a finite number of them. Note that letting $\alpha=0$ we obtain that the line in \eqref{1-20} becomes the real axis, since $\nu_-=-\nu_+$ and hence $m_-=m_+$ and $\rho_-=\rho_+$ in this case. Of course, this is consistent with the trivial zeros of $F(s)$, which lie on the real axis thanks to the special form of \eqref{1-6}, and with asymptotic formula \eqref{1-21}.

\medskip
{\bf Remark 4.} Note that when $\alpha$ varies in such a way that $\sqrt{n}<\nu_\alpha<\sqrt{n+1}$ with some $n\in\NN$, the coefficients of the line in \eqref{1-20} may in principle continuously vary between $-\infty$ and $+\infty$. But when $\nu_\alpha$ hits the square root of an integer, the position of such a line may change suddenly. This shows an interesting discontinuity in the $\alpha$-behavior of $F(s,\alpha)$.  \fine

\medskip
{\bf 1.7.~Nontrivial zeros.} Now we turn to the nontrivial zeros of $F(s,\alpha)$. Let $\si^+$ be the upper bound of the real parts of the zeros of $F(s,\alpha)$ and $\si^-$ be as in the previous subsection. The nontrivial zeros of $F(s,\alpha)$ are those in the vertical strip $-\si^-\leq \si \leq \si^+$, thus the zeros of $F(s,\alpha)$ are the disjoint union of trivial and nontrivial zeros. Let
\begin{equation}
\label{1-22}
N_F(T,\alpha) = |\{\rho=\beta+i\gamma: F(\rho,\alpha)=0, \si^-\leq \beta \leq \si^+, |\gamma|\leq T\}|
\end{equation}
be the counting function of the nontrivial zeros, and let $\overline{n}$ be the smallest $n$ with $a(n)\neq0$. A suitable application of the argument principle gives the following analog of the Riemann-von Mangoldt formula.

\medskip
{\bf Corollary 3.} {\sl Let $\alpha>0$. Then as $T\to\infty$ we have}
\begin{equation}
\label{1-23}
N_F(T,\alpha) = \frac{2}{\pi}T\log T \ + \ \frac{T}{\pi} \log\big(\frac{N}{\overline{n}m_+m_-(2\pi e)^2}\big) \ + \ O(\log T).
\end{equation}

\medskip
Again, note that letting $\alpha=0$ in Corollary 3 we get the well known asymptotic formula for the counting function of the nontrivial zeros of $F(s)$, obtained substituting $m_+m_-$ by $\overline{n}$ in \eqref{1-23}, since indeed $m_+m_- = \overline{n}$ in this case. Also, thanks to the constant 2 in the main term of \eqref{1-23}, Corollary 3 shows once more the degree 2 behavior of $F(s,\alpha)$. Further, when $\alpha$ varies we can observe a similar phenomenon as in Remark 4 in the behavior of the coefficient of $T$ in \eqref{1-23}.

\medskip
We conclude remarking that several other problems can be studied once the functional equation is established, for example Voronoi type formulae and related estimates for nonlinear exponential sums, with modular coefficients in the case treated in this paper. We shall return on these and other questions in a future paper. Anyway, already the above theorem shows that the standard twist $F(s,\alpha)$ is a respectable object from the $L$-functions point of view, at least in the case of half-integral weight modular forms. Moreover, although the arguments in this paper are not sufficient to prove the functional equation of $F(s,\alpha)$ for general $F\in\S^\sharp$, we believe that indeed $F(s,\alpha)$ satisfies a suitable functional equation in such a general case. Again, we shall return on this question in a future paper, but here we add a final remark.

\medskip
{\bf Remark 5.} Thanks to the characterization in \cite{Ka-Pe/1999a} of the degree 1 functions of $\S^\sharp$ as linear combinations of Dirichlet $L$-functions over Dirichlet polynomials from $\S^\sharp$, the standard twists of degree 1 functions in $\S^\sharp$ are closely related to the Hurwitz-Lerch zeta function
\[
L(s,\alpha,\lambda) = \sum_{n=0}^\infty \frac{e(-n\alpha)}{(n+\lambda)^s} \hskip1.5cm 0\leq \alpha<2\pi, \ 0<\lambda\leq1.
\]
Since such zeta functions satisfy a functional equation, see e.g. Berndt \cite{Ber/1972}, one certainly expects that the standard twits of degree 1 functions satisfy a functional equation of Hurwitz-Lerch type. Since \eqref{1-16} may clearly be regarded as a degree 2 analog of the Hurwitz-Lerch functional equation, and the same holds in the higher degree cases mentioned in Remark 2, it is not unreasonable to expect that, in general, standard twists satisfy a kind of general functional equation of Hurwitz-Lerch type. \fine

\medskip
{\bf Acknowledgements.} 
We wish to thank Jan Hendrik Bruinier and Winfried Kohnen for their advice concerning half-integral weight modular forms. 
This research was partially supported by the Istituto Nazionale di Alta Matematica, by the MIUR grant PRIN-2015 {\sl ``Number Theory and Arithmetic Geometry''} and by grant 2017/25/B/ST1/00208 {\sl ``Analytic methods in number theory''}  from the National Science Centre, Poland.

\bigskip
\section{Proof of the Theorem and Corollary 1}

\smallskip
{\bf 2.1.~Basic formula.} For convenience we write 
\begin{equation}
\label{2.1}
Q = \sqrt{N}/2\pi,
\end{equation}
thus, recalling \eqref{1-10}, functional equation \eqref{1-6} becomes
\begin{equation}
\label{2.2}
Q^s\Gamma(s+\mu)F(s) = \omega Q^{1-s}\Gamma(1-s+\mu)F^*(1-s).
\end{equation}
By the reflection formula $\Gamma(z)^{-1} =  \Gamma(1-z) \sin(\pi z)/\pi$ with $z=s+\mu$, we transform \eqref{2.2} in the asymmetric form
\[
F(s) = \frac{\omega Q^{1-2s}}{\pi} \Gamma(1-s+\mu)\Gamma(1-s-\mu) \sin\pi(s+\mu) F^*(1-s).
\]
Therefore, writing for $h^*$ as in \eqref{1-10}
\begin{equation}
\label{2.300}
\mu^*= \frac{2h^*+1}{4} 
\end{equation}
and observing that $\mu=\pm\mu^*$ if $h^*=0$ and $\mu=\mu^*$ if $h^*\geq1$, for every $h^*\geq0$  we rewrite the above functional equation as
\begin{equation}
\label{2.3}
F(s) = \frac{\omega Q^{1-2s}}{\pi} \Gamma(1-s+\mu^*)\Gamma(1-s-\mu^*) \sin\pi(s+\mu) F^*(1-s).
\end{equation}
Thanks to the special value of $\mu$ (and hence of $\mu^*$) in \eqref{1-10} we may suitably transform the two $\Gamma$-factors in \eqref{2.3}. To this and, since $\mu^*-h^*=-\mu^*+1/2$, we first apply $h^*$-times the factorial formula to the first $\Gamma$-factor in \eqref{2.3} to get (recall that an empty product equals 1)
\[
\Gamma(1-s+\mu^*) = \Gamma(1-s-\mu^*+1/2) \prod_{1\leq j \leq h^*}(1-s+\mu^*-j).
\]
Then we apply the duplication formula $\Gamma(z) \Gamma(z+1/2) = \pi^{1/2} 2^{1-2z} \Gamma(2z)$ with $z=1-s-\mu^*$. In view of \eqref{2.300}, a simple computation shows that \eqref{2.3} becomes
\begin{equation}
\label{2.4}
\begin{split}
F(s) = \frac{2\omega}{\sqrt{2\pi}} &\big(\frac{Q}{2}\big)^{1-2s} \Gamma(2(1-s)-1/2-h^*) \\
&\times \left(\prod_{1\leq j \leq h^*} (2(1-s)-1/2-h^*+2j-1)\right) \sin\pi(s+\mu) F^*(1-s).
\end{split}
\end{equation}
Applying identity \eqref{1-15} with $X=2(1-s)-1/2-h^*$ we see that \eqref{2.4} can be written as
\begin{equation}
\label{2.5}
F(s) = \frac{2\omega}{\sqrt{2\pi}} \big(\frac{Q}{2}\big)^{1-2s} \left(\sum_{\ell=0}^{h^*} a_\ell \Gamma(2(1-s)-1/2-\ell)\right) \sin\pi(s+\mu) F^*(1-s).
\end{equation}

\smallskip
Let now 
\begin{equation}
\label{2.6}
X > 1, \qquad  \alpha>0, \qquad z_X=\frac{1}{X}+2\pi i\alpha
\end{equation}
and
\begin{equation}
\label{2.7}
F_X(s,\alpha) = \sum_{n=1}^\infty \frac{a(n)}{n^s}e(-\alpha \sqrt{n})e^{-\sqrt{n}/X} =  \sum_{n=1}^\infty \frac{a(n)}{n^s} e^{-z_X \sqrt{n}};
\end{equation}
clearly, $F_X(s,\alpha)$ converges for every $s\in\CC$.
Given $K>0$, for $-K<\si<0$ by Mellin's transform we have
\[
F_X(s,\alpha) = \frac{1}{2\pi i}\int_{(2(K+2))} F(s+\frac{w}{2}) \Gamma(w) z_X^{-w} \d w,
\]
then we shift the line of integration to $\Re w=\delta$; here and later, $\delta>0$ is a sufficiently small constant, not necessarily the same at each occurrence. Since $F(s)$ is entire, for $-K<\sigma<-\delta$ we may use functional equation \eqref{2.5} and expand $F^*(1-s-w/2)$ to obtain
\begin{equation}
\label{2.8}
\begin{split}
F_X(s,\alpha) &= \frac{2\omega}{\sqrt{2\pi}} \big(\frac{Q}{2}\big)^{1-2s} \sum_{n=1}^\infty \frac{a^*(n)}{n^{1-s}} \sum_{\ell=0}^{h^*} a_\ell \\
&\times \frac{1}{2\pi i} \int_{(\delta)} \Gamma(2(1-s)-1/2-\ell-w) \Gamma(w) \sin\pi(s+\frac{w}{2}+\mu) \big(\frac{Qz_X}{2\sqrt{n}}\big)^{-w} \d w.
\end{split}
\end{equation}
Using the expression $\sin \pi z = (e^{i\pi z}-e^{-i\pi z})/2i$ we rewrite \eqref{2.8} as
\begin{equation}
\label{2.9}
\begin{split}
F_X(s,\alpha) &= \frac{\omega}{i\sqrt{2\pi}}  \big(\frac{Q}{2}\big)^{1-2s} \sum_{n=1}^\infty \frac{a^*(n)}{n^{1-s}}  \sum_{\ell=0}^{h^*} a_\ell \\
&\times \big(e^{i\pi(s+\mu)} H(s+\ell/2,z_X(n,\alpha)) - e^{-i\pi(s+\mu)} H(s+\ell/2,-z_X(n,\alpha))\big),
\end{split}
\end{equation}
where
\begin{equation}
\label{2.10}
z_X(n,\alpha) = \frac{\pi \alpha Q}{\sqrt{n}} -i\frac{Q}{2X\sqrt{n}}
\end{equation}
in view of \eqref{2.6}, and
\begin{equation}
\label{2.11}
H(s,z) = \frac{1}{2\pi i} \int_{(\delta)} \Gamma(2(1-s)-1/2-w) \Gamma(w) z^{-w} \d w.
\end{equation}

\smallskip
An explicit expression for the function $H(s,z)$ can be obtained from the formula for the Mellin-Barnes integral in (3.3.9) of Ch.3 of Paris-Kaminski \cite{Pa-Ka/2001}, namely
\begin{equation}
\label{2.12}
\frac{1}{2\pi i} \int_{(c)} \Gamma(\xi-w) \Gamma(w) z^{-w} \d w = \Gamma(\xi) (1+z)^{-\xi}
\end{equation}
provided $0<c<\Re \xi$ and $|\arg z |<\pi$. To evaluate $H(s+\ell/2,\pm z_X(n,\alpha))$ we have to check the first condition, the second one being certainly satisfied in our case. Recalling that we already have the condition $-K<\si<-\delta$, we choose $K=h^*+10$ and
\begin{equation}
\label{2.13}
-2\delta<\si<-\delta \ \ \text{if} \ \ h^* \leq 2 \ \ \text{and} \ \ s_{h^*+1}+\delta <\si<s_{h^*}-\delta \ \ \text{if} \ \ h^*\geq 3,
\end{equation}
so in view of \eqref{1-8} the first condition required by \eqref{2.12} is also satisfied for every $0\leq \ell \leq h^*$. Therefore, from \eqref{2.11} and \eqref{2.12} we get
\begin{equation}
\label{2.14}
H(s+\ell/2,\pm z_X(n,\alpha)) = \Gamma(2(1-s)-1/2-\ell) (1\pm z_X(n,\alpha))^{2(s-s_\ell)},
\end{equation}
and inserting \eqref{2.14} into \eqref{2.9} we obtain, for $\sigma$ as in \eqref{2.13}, the basic formula
\begin{equation}
\label{2.15}
\begin{split}
F_X(s,\alpha) &=  \frac{\omega}{i\sqrt{2\pi}} \big(\frac{Q}{2}\big)^{1-2s} \sum_{\ell=0}^{h^*} a_\ell  \Gamma(2(1-s)-1/2-\ell)  \sum_{n=1}^\infty \frac{a^*(n)}{n^{1-s}} \times \\
&\times \big\{e^{i\pi(s+\mu)} (1+ z_X(n,\alpha))^{2(s-s_\ell)} - e^{-i\pi(s+\mu)} (1-z_X(n,\alpha))^{2(s-s_\ell)}\big\},
\end{split}
\end{equation}
where $\omega$, $a_\ell$, $z_X(n,\alpha)$ and $s_\ell$ are given by \eqref{1-6}, \eqref{1-15}, \eqref{2.10} and \eqref{1-8}, respectively. Note that the series in \eqref{2.15} is absolutely convergent since $\si<0$ and $|(1\pm z_X(n,\alpha))^{2(s-s_\ell})|\to 1$ as $n\to \infty$.

\bigskip
{\bf 2.2.~Limit as $X\to\infty$.} Letting $X\to\infty$ in \eqref{2.15} requires some care. Indeed, the limit of $F_X(s,\alpha)$ is clearly $F(s,\alpha)$ for every $\alpha>0$ when $\si>1$, but \eqref{2.15} holds in a different range, namely for $s$ as in  \eqref{2.13}. Moreover, the limit of the term $(1- z_X(n,\alpha))^{2(s-s_\ell)}$ is not always well defined, since as $X\to\infty$ we have $1- z_X(n,\alpha)\to 1- \pi\alpha Q/\sqrt{n}$, which may vanish for suitable values of $\alpha$ and $n$. Actually, in view of \eqref{2.1} and of the definition of $n_\alpha$ and Spec$^*(F)$ in \eqref{1-7}, such critical values of $\alpha$ and $n$ arise when $\alpha\in$ Spec$^*(F)$ and $n=n_\alpha$. Therefore, before letting $X\to\infty$ we derive a different expression for $F_X(s,\alpha)$.

\smallskip
Since we already have information on the analytic properties of $F(s,\alpha)$, see Subsection 1.1 or \cite{Ka-Pe/2005} and \cite{Ka-Pe/resoI}, we now consider $F_X(s,\alpha)$ as the twist of $F(s,\alpha)$ by $e^{-\sqrt{n}/X}$, see \eqref{2.7}. Hence, by Mellin's transform, for $s$ as in \eqref{2.13}, $X>1$ and $K=h^*+10$ we have
\[
F_X(s,\alpha) = \frac{1}{2\pi i}\int_{(2(K+2))} F(s+\frac{w}{2},\alpha) \Gamma(w) X^w \d w.
\]
The integrand has simple poles at $w=-r$, $r=0,1,\dots$ and, as already mentioned in the Introduction, at most simple poles at $w=2(s_\ell-s)$, where $s_\ell$ are defined by \eqref{1-8}. We denote by $\rho_\ell(\alpha)$ the residue of $F(s+w/2,\alpha)$ at $w=2(s_\ell-s)$; clearly
\begin{equation}
\label{2.151}
\rho_\ell(\alpha) =2\kappa_\ell(\alpha)
\end{equation}
where $\kappa_\ell(\alpha)$ is the residue of $F(s,\alpha)$ at $s_\ell$. We shift the line of integration to $\Re w=-\delta$ (once more, with a sufficiently small $\delta>0$, not necessarily the same as at previous places), thus collecting residues of the poles at $w=0$ and, recalling that $s$ is as in \eqref{2.13}, also at $w= 2(s_\ell-s)$ with $0\leq \ell\leq h^*$. Therefore, still for $s$ as in \eqref{2.13}, we get
\begin{equation}
\label{2.16}
\begin{split}
F_X(s,\alpha) &= F(s,\alpha) + \sum_{\ell=0}^{h^*} \rho_\ell(\alpha) \Gamma(2(s_\ell-s)) X^{2(s_\ell-s)} \\
&+ \frac{1}{2\pi i}\int_{(-\delta)} F(s+\frac{w}{2},\alpha) \Gamma(w) X^w \d w \\
&=  F(s,\alpha) + \Sigma_X(s,\alpha) + I_X(s,\alpha),
\end{split}
\end{equation}
say. Moreover, since $F(s,\alpha)$ has polynomial growth, see \cite{Ka-Pe/2005}, \cite{Ka-Pe/resoI} or the Introduction, as $X\to\infty$  we have
\begin{equation}
\label{2.17}
I_X(s,\alpha) \ll X^{-\delta} \int_{(-\delta)} |F(s+\frac{w}{2},\alpha) \Gamma(w) \d w| \to 0.
\end{equation}

\smallskip
On the other hand, rewriting \eqref{2.15} as (recall that $a^*(n_\alpha)=0$ if $n_\alpha\not\in\NN$)
\begin{equation}
\label{2.18}
\begin{split}
F_X(s,\alpha) &=  \frac{\omega}{i\sqrt{2\pi}} \big(\frac{Q}{2}\big)^{1-2s} \sum_{\ell=0}^{h^*} a_\ell  \Gamma(2(1-s)-1/2-\ell) \times \\
&\hskip2cm \times  \big\{e^{i\pi(s+\mu)} \sum_{n=1}^\infty \frac{a^*(n)}{n^{1-s}} (1+ z_X(n,\alpha))^{2(s-s_\ell)} \\
&\hskip2cm -  e^{-i\pi(s+\mu)} \sum_{\substack{n=1 \\ n\neq n_\alpha}}^\infty \frac{a^*(n)}{n^{1-s}}  (1-z_X(n,\alpha))^{2(s-s_\ell)} \\
&\hskip2cm -  e^{-i\pi(s+\mu)}  \frac{a^*(n_\alpha)}{n_\alpha^{1-s}}  (1-z_X(n_\alpha,\alpha))^{2(s-s_\ell)}\big\} \\
&= A_X(s,\alpha) - B_X(s,\alpha) - C_X(s,\alpha)
\end{split}
\end{equation}
say, from \eqref{2.16} and \eqref{2.18} we get, for $X>1$ and $s$ as in \eqref{2.13}, that
\begin{equation}
\label{2.19}
F(s,\alpha)  = A_X(s,\alpha) - B_X(s,\alpha) - \big(C_X(s,\alpha) + \Sigma_X(s,\alpha)\big) - I_X(s,\alpha).
\end{equation}
From \eqref{2.17} and \eqref{2.18} it is clear that, as $X\to\infty$, $I_X(s,\alpha)$, $A_X(s,\alpha)$ and $B_X(s,\alpha)$ tend, respectively, to $0$ and to well defined functions $A(s,\alpha)$ and $B(s,\alpha)$, say, in the range \eqref{2.13}. Later on we shall compute explicitly such functions, but first we treat the critical term $C_X(s,\alpha)+\Sigma_X(s,\alpha)$, depending on $\alpha$ belonging to Spec$^*(F)$ or not. Actually, we show that in both cases $C_X(s,\alpha)+\Sigma_X(s,\alpha)$ vanishes identically.

\smallskip
Suppose first that $\alpha\not\in$ Spec$^*(F)$. Then $F(s,\alpha)$ is entire and $a^*(n_\alpha)=0$, hence both $C_X(s,\alpha)$ and $\Sigma_X(s,\alpha)$ vanish identically. Suppose now $\alpha\in$ Spec$^*(F)$. Then from \eqref{1-7}, \eqref{2.1} and \eqref{2.10} we have
\[
1- z_X(n_\alpha,\alpha) = i\frac{Q}{2X\sqrt{n_\alpha}}.
\]
Since in view of \eqref{1-8} we have $2(s_\ell-s) = 2(1-s)-1/2-\ell$, and recalling the definition of conductor $q$ of $F(s)$ in \eqref{1-3}, by a simple computation we obtain
\[
C_X(s,\alpha) = \sum_{\ell=0}^{h^*} \pi_\ell(\alpha) \Gamma(2(s_\ell-s)) X^{2(s_\ell - s)} 
\]
with
\begin{equation}
\label{2.20}
\pi_\ell(\alpha) = -\frac{i\omega a_\ell}{\sqrt{2\pi}} \frac{a^*(n_\alpha)}{n_\alpha^{1-s_\ell}}\big(\frac{\sqrt{N}}{4\pi}\big)^{1-2s_\ell} e^{-i\pi(s_\ell+\mu)}.
\end{equation}
Since $\Re(s_\ell-s)>0$ for $s$ as in \eqref{2.13}, comparing $\Sigma_X(s,\alpha)$ in \eqref{2.16} with the above expression for $C_X(s,\alpha)$ and observing that all other terms in \eqref{2.19} remain bounded as $X\to\infty$, we deduce that $C_X(s,\alpha)+\Sigma_X(s,\alpha)$ must vanish identically. In particular,  thanks to \eqref{2.151} and \eqref{2.20} the residue $\kappa_\ell(\alpha)$ of $F(s,\alpha)$ at $s=s_\ell$ is explicitly given by
\begin{equation}
\label{2.21}
\kappa_\ell(\alpha) =  \frac{i\omega a_\ell}{2\sqrt{2\pi}} \frac{a^*(n_\alpha)}{n_\alpha^{1-s_\ell}}\big(\frac{\sqrt{N}}{4\pi}\big)^{1-2s_\ell} e^{-i\pi(s_\ell+\mu)} \hskip1.5cm \ell=0,\dots,h^*,
\end{equation}
where $\omega$, $a_\ell$ and $\mu$ are defined by \eqref{1-4}, \eqref{1-15} and \eqref{1-10}, respectively.

\smallskip
Finally we may let $X\to\infty$ in \eqref{2.19} thus getting, for every $\alpha>0$ and $s$ as in \eqref{2.13}, that
\begin{equation}
\label{2.22}
F(s,\alpha) = A(s,\alpha) - B(s,\alpha).
\end{equation}

\bigskip
{\bf 2.3.~Functional equation.} We first compute $A(s,\alpha)$ and $B(s,\alpha)$. Note that by  \eqref{1-7}, \eqref{1-11} and \eqref{2.1}
\[
\pi\alpha Q = \sqrt{n_\alpha} = \nu_\alpha.
\] 
Hence, letting $X\to\infty$, from \eqref{2.10} and  \eqref{2.18} we obtain at once that
\begin{equation}
\label{2.23}
A(s,\alpha) = \frac{\omega}{i\sqrt{2\pi}} \big(\frac{Q}{2}\big)^{1-2s} \sum_{\ell=0}^{h^*} a_\ell  \Gamma(2(1-s)-1/2-\ell) e^{i\pi(s+\mu)} \sum_{n=1}^\infty \frac{a^*(n)}{n^{1-s}} (1+ \frac{\nu_\alpha}{\sqrt{n}})^{2(s-s_\ell)}.
\end{equation}
From \eqref{1-7} and \eqref{2.10} we see that for $n<n_\alpha$ the real part of $1-z_X(n,\alpha)$ is negative, hence for such $n$'s we have
\[
\lim_{X\to\infty} (1- z_X(n,\alpha))^{2(s-s_\ell)} = \left|1-\frac{\nu_\alpha}{\sqrt{n}}\right|^{2(s-s_\ell)} e^{i\pi 2(s-s_\ell)}.
\]
Therefore, arguing as for $A(s,\alpha)$ and recalling that an empty sum equals 0, we get
\begin{equation}
\label{2.24}
\begin{split}
B(s,\alpha) &= \frac{\omega}{i\sqrt{2\pi}} \big(\frac{Q}{2}\big)^{1-2s}  \sum_{\ell=0}^{h^*} a_\ell  \Gamma(2(1-s)-1/2-\ell) e^{-i\pi(s+\mu)} \times\\
& \times\left\{ e^{i\pi 2(s-s_\ell)} \sum_{1\leq n<n_\alpha} \frac{a^*(n)}{n^{1-s}} \left|1- \frac{\nu_\alpha}{\sqrt{n}}\right|^{2(s-s_\ell)}
+ \sum_{n>n_\alpha} \frac{a^*(n)}{n^{1-s}} \left|1- \frac{\nu_\alpha}{\sqrt{n}}\right|^{2(s-s_\ell)}\right\}.
\end{split}
\end{equation}
Note, as after \eqref{2.15}, that the series in \eqref{2.23} and \eqref{2.24} are absolutely convergent since $\si<0$.

\smallskip
Next we give a uniform shape to the generalized Dirichlet series in $A(s,\alpha)$ and $B(s,\alpha)$. To accomplish this we recall that $\nu=\pm\sqrt{n}$ with $n=1,2\dots$, hence for $\nu\neq -\nu_\alpha$ we have
\begin{equation}
\label{2.25}
\frac{1}{|\nu|^{2(1-s_\ell)} |\nu+\nu_\alpha|^{2(s_\ell-s)}} =
\begin{cases}
\frac{1}{n^{1-s}} \big(1+\frac{\nu_\alpha}{\sqrt{n}}\big)^{2(s-s_\ell)} & \text{if} \ \nu\geq 1 \\
\frac{1}{n^{1-s}} \left|1-\frac{\nu_\alpha}{\sqrt{n}}\right|^{2(s-s_\ell)} & \text{if} \ \nu\leq -1.
\end{cases}
\end{equation}
Therefore, a simple computation based on \eqref{2.22}-\eqref{2.25} and \eqref{1-8} shows that for $s$ as in \eqref{2.13} we have
\begin{equation}
\label{2.28}
F(s,\alpha) = \frac{\omega}{i\sqrt{2\pi}} \big(\frac{Q}{2}\big)^{1-2s} \sum_{\ell=0}^{h^*} a_\ell  \Gamma(2(1-s)-1/2-\ell) F^*_\ell(1-s,\alpha)
\end{equation}
with $F^*_\ell(s,\alpha)$ as in \eqref{1-12}. The Theorem now follows by analytic continuation from \eqref{2.1}, \eqref{2.28} and the following lemma.

\medskip
{\bf Lemma 2.1.} {\sl Let $\alpha>0$ and $\ell=0,\dots,h^*$. Then the functions $F_\ell^\pm(s,\alpha)$ are entire with polynomial growth on vertical strips.} 

\medskip
{\it Proof.} We prove the lemma by showing that the $F_\ell^\pm(s,\alpha)$'s are close to certain stratifications of $F^*(s)$. Let  $H>\nu_\alpha^2+1$ be an integer to be chosen later on. From \eqref{2.25} and the definition of $F^\pm_\ell(s,\alpha)$, see \eqref{1-14}, for $\si>1$ we have
\begin{equation}
\label{2.29}
F^\pm_\ell(s,\alpha) = \sum_{n\geq H} \frac{c^*(n)}{n^s} \big(1\pm\frac{\nu_\alpha}{\sqrt{n}}\big)^{-2s+1/2+\ell} + G_\ell^\pm(s)
\end{equation}
with certain generalized Dirichlet polynomials $G_\ell^\pm(s)= G^\pm_{\ell,H}(s,\alpha)$, which are entire. Writing 
\begin{equation}
\label{2.290}
\rho=-2s+1/2+\ell
\end{equation}
we have
\[
\big(1\pm \frac{\nu_\alpha}{\sqrt{n}}\big)^\rho= \sum_{r=0}^\infty (\pm1)^r{\rho\choose r}\big(\frac{\nu_\alpha}{\sqrt{n}}\big)^r = \Sigma_1 + \Sigma_2,
\]
where $\Sigma_1$ is the sum over $r\leq R$ and $\Sigma_2$ over $r\geq R+1$, $R$ being an arbitrarily large positive integer. Clearly
\[
\big|(\pm1)^r{\rho \choose r}\big| \leq \frac{|\rho||\rho+1|\cdots|\rho+r-1|}{r!} \leq \prod_{j=1}^r(1+|\rho|/j) \leq (1+|\rho|)^r,
\]
therefore, writing 
\begin{equation}
\label{2.291}
K=K_{s,\ell,\alpha} = (1+|\rho|)\nu_\alpha
\end{equation}
and choosing $H\geq2K^2$, for $n\geq H$ we obtain
\[
\Sigma_2  \ll  \sum_{r=R+1}^\infty (\frac{K}{\sqrt{n}})^r = (\frac{K}{\sqrt{n}})^{R+1} \frac{1}{1-K/\sqrt{n}} \ll (\frac{K}{\sqrt{n}})^{R+1}.
\]
Consequently we have
\begin{equation}
\label{2-37}
\begin{split}
\sum_{n\geq H} \frac{c^*(n)}{n^s} \big(1\pm\frac{\nu_\alpha}{\sqrt{n}}\big)^\rho &= 
 \sum_{n\geq H} \frac{c^*(n)}{n^s} \sum_{r\leq R}(\pm1)^r{\rho\choose r} \big(\frac{\nu_\alpha}{\sqrt{n}}\big)^r \\
&+O\big(K^{R+1}\sum_{n\geq H} \frac{|c^*(n)|}{n^{\sigma+(R+1)/2}} \big) = E^\pm_1(s) + E^\pm_2(s),
\end{split}
\end{equation}
say, where $E^\pm_j(s) = E^\pm_{j,\ell,H,R}(s,\alpha)$. 

\medskip
Let $\K$ be any compact subset of $\CC$ intersecting the half-plane $\si>1$ and choose
\begin{equation}
\label{2.292}
H= H_{\K,\ell,\alpha} = 2[\max_{s\in\K} K^2]+1.
\end{equation}
Then $E^\pm_2(s)$ is a double series of holomorphic functions, absolutely and uniformly convergent in $\{\si> -(R-1)/2\}\cap\K$. Now observe that, in view of \eqref{1-13}, for $n\geq H$ we have, depending on the the choice of $\pm$ in $F_\ell^\pm(s,\alpha)$, 
\[
c^*(n) = \mp e^{\pm i\pi\mu} a^*(n).
\]
Therefore, rearranging terms we obtain
\begin{equation}
\label{2-38}
\begin{split}
E^\pm_1(s) &= \mp e^{\pm i\pi\mu} \sum_{r\leq R} (\pm1)^r {\rho\choose r} \nu_\alpha^r \sum_{n\geq H} \frac{a^*(n)}{n^{s+r/2}}
\\
&= \mp e^{\pm i\pi\mu} \sum_{r\leq R} (\pm1)^r {\rho\choose r} \nu_\alpha^r \big(F^*(s+r/2) - D^*(s+r/2)\big) \\
&=\sum_{r\leq R} Q^\pm_r(s) F^*(s+r/2) + E^\pm_3(s),
\end{split}
\end{equation}
where $D^*(s)=D^*_H(s)$ is a Dirichlet polynomial,  the $Q^\pm_r(s)= Q^\pm_{r,\ell}(s,\alpha)$ are polynomials and $E^\pm_3(s)=E^\pm_{3,\ell,H,R}(s,\alpha)$ is an entire function. Thus $E^\pm_1(s)$ is also entire thanks to the properties of $F^*(s)$, and hence from \eqref{2.29}-\eqref{2-38} we deduce that the generalized Dirichlet series $F^\pm_\ell(s,\alpha)$, absolutely convergent for $\si>1$, are actually entire functions since $R$ and $\K$ are arbitrary.

\smallskip
To prove that $F^\pm_\ell(s,\alpha)$ have polynomial growth on vertical strips we write, in view of \eqref{2.29}, \eqref{2.290} and \eqref{2-37},
\begin{equation}
\label{2-39}
F^\pm_\ell(s,\alpha) = G_\ell^\pm(s) + E^\pm_1(s) + E^\pm_2(s)
\end{equation}
and recall that the functions on the right hand side depend on $H$. Moreover, from \eqref{2.290}, \eqref{2.291}, \eqref{2.292} and choosing $\K$ to be the rectangle $[-(R-2)/2 , 2]\times [-|t|, |t|]$, we have that
\[
H \ll_R (|t|+1)^2
\]
uniformly for $-(R-2)/2\leq \si\leq 2$. Since the (generalized) Dirichlet polynomials involved in \eqref{2-39} have length $\ll H$ and coefficients closely related with $a^*(n)$, recalling that $F^*(s)$ has polynomial growth on vertical strips we see
from \eqref{2-37}, \eqref{2-38} and \eqref{2-39} that
\[
F^\pm_\ell(s,\alpha) \ll_R (|t|+1)^{c(R)}
\]
with some $c(R)>0$, uniformly for $-(R-2)/2\leq \si\leq 2$ with an arbitrarily large $R>0$. \fine

\bigskip
{\bf 2.4.~Proof of Corollary 1.} The first assertion is already proved in Lemma 2.1. We note preliminarily that from \eqref{1-12} we have
\[
\inf\{\xi: e^{-\pi |t|} |F_\ell^*(\sigma+it,\alpha)| \ll |t|^\xi \ \text{as} \ |t|\to\infty\} \geq 0
\]
for every $\si$, since $F_\ell^\pm(s,\alpha)$ are generalized Dirichlet series with polynomial growth on vertical strips. Hence Stirling's formula coupled with \eqref{1-16} shows that 
\begin{equation}
\label{2-43}
\mu(\si) \geq 1-2\si,
\end{equation}
since the factorial formula of the $\Gamma$ function implies that for $\ell=0,\dots,h^*$
\begin{equation}
\label{2-44}
\Gamma(2(1-s)-1/2-\ell) = \Gamma(2(1-s)-1/2)/P_\ell(s)
\end{equation}
with certain polynomials $P_\ell(s)$ of degree $\ell$ and $P_0(s)\equiv 1$. Now we proceed in a similar way as in Lemma 2.1 to show that $F_\ell^\pm(s,\alpha)$ are a kind of stratification of $F_0^\pm(s,\alpha)$. Clearly we may write $|\nu||\nu+\nu_\alpha|$, with $\nu=\pm\sqrt{n}$ and $n$ large enough, as $\nu(\nu\pm\nu_\alpha)$ with $\nu=\sqrt{n}$. Hence, choosing  $V>0$ sufficiently large, for $\ell=0,\dots,h^*$ we write
\begin{equation}
\label{2-40}
F_\ell^{\pm}(s,\alpha)= \sum_{\nu\geq V}\frac{c^*(\nu^2)}{\nu^{1/2+\ell}(\nu\pm\nu_\alpha)^{2s-1/2-\ell}} + D_\ell^{\pm}(s,\alpha) = G_\ell^{\pm}(s,\alpha) + D_\ell^{\pm}(s,\alpha),
\end{equation}
say, where $D_\ell^{\pm}(s,\alpha)$ are generalized Dirichlet polynomials. Moreover, since
\[
\left(\frac{\nu}{\nu\pm\nu_\alpha}\right)^{-\ell} = \left(1 \mp \frac{\nu_\alpha}{\nu\pm\nu_\alpha}\right)^{-\ell} = \sum_{r=0}^\infty (\mp\nu_\alpha)^r {-\ell \choose  r} \frac{1}{(\nu\pm\nu_\alpha)^r},
\]
we have
\[
G_\ell^{\pm}(s,\alpha) = \sum_{r=0}^\infty (\mp\nu_\alpha)^r {-\ell \choose  r} G_0^{\pm}(s+r/2,\alpha)
\]
and therefore by \eqref{2-40}
\[
F_\ell^{\pm}(s,\alpha) = \sum_{r=0}^\infty (\mp\nu_\alpha)^r {-\ell \choose  r} \left(F_0^{\pm}(s+r/2,\alpha) - D_0^{\pm}(s+r/2,\alpha)\right) + D_\ell^{\pm}(s,\alpha).
\]
Hence, since the generalized Dirichlet series are bounded in the half-plane of absolute convergence, for every $\sigma$ there exists $K=K(\sigma)$ such that as $|t|\to\infty$
\begin{equation}
\label{2-41}
F_\ell^{\pm}(1-s,\alpha) = \sum_{r=0}^K(\mp\nu_\alpha)^r {-\ell \choose  r} F_0^{\pm}(1-s+r/2,\alpha) + O(1).
\end{equation}

\smallskip
From \eqref{1-12}, \eqref{1-16} and \eqref{2-41} we therefore have
\begin{equation}
\label{2-42}
F(s,\alpha) = e^{as+b}\sum_{\ell=0}^{h^*} a_\ell \Gamma(2(1-s)-1/2-\ell)\left\{H_\ell(1-s,\alpha) +O(e^{\pi |t|})\right\}
\end{equation}
with suitable $a\in\RR$, $b\in\CC$ and
\[
\begin{split}
H_\ell(1-s,\alpha) &= e^{-i\pi(1-s)}  \sum_{r=0}^K(-\nu_\alpha)^r {-\ell \choose  r} F_0^+(1-s+r/2,\alpha) \\
&+ e^{i\pi(1-s)}  \sum_{r=0}^K(\nu_\alpha)^r {-\ell \choose  r} F_0^-(1-s+r/2,\alpha) \\
& = H^+_\ell(1-s,\alpha)  + H^-_\ell(1-s,\alpha),
\end{split}
\]
say. Note that if $t\to+\infty$, then $H^+_\ell(1-s,\alpha)$ has exponential decay since $F_0^+(s,\alpha)$ has polynomial growth, and similarly for $H^-_\ell(1-s,\alpha)$ if $t\to-\infty$. Suppose that $\mu^*(1-\si)=0$, and hence by \eqref{1-18} also $\mu^\pm(1-\si)=0$. Then, again by Stirling's formula and \eqref{2-44}, from \eqref{2-42} we get that $\mu(\si) \leq 1-2\si$ and therefore by \eqref{2-43}
\[
\mu(\si) = 1-2\si = 1-2\si +\mu^*(1-\si)
\]
in this case. Suppose now that  $\mu^*(1-\si)>0$. Then, thanks to definitions \eqref{1-17} and \eqref{1-18}, the error term in \eqref{2-42} is negligible and, still by Stirling's formula and \eqref{2-44}, from \eqref{2-42} we get
\[
\inf\{\xi:F(s,\alpha)\ll |t|^\xi\} = 
\begin{cases}
1-2\si +\mu^-(1-\si) \ \  \text{if} \ \ t\to+\infty\\
1-2\si +\mu^+(1-\si) \ \  \text{if} \ \ t\to-\infty.
\end{cases}
\]
The proof of the corollary is now complete. \fine

\bigskip
\section{Proof of Corollaries 2 and 3}

\smallskip
{\bf 3.1.~Proof of Corollary 2.} We start recalling the notation of Subsection 1.6 and writing
\[
\Sigma(s) = \sum_{\ell=0}^{h^*}a_\ell \Gamma\big(2s-1/2-\ell\big)F^*_\ell(s,\alpha);
\]
note that by \eqref{1-16} the zeros of $F(s,\alpha)$ coincide with those of $\Sigma(1-s)$. By \eqref{2-44} we write, again with obvious polynomials $Q_\ell(s)$ of degree $\ell$,
\begin{equation}
\label{3-1}
\Sigma(s) = \Gamma(2s-1/2) \sum_{\ell=0}^{h^*}a_\ell \frac{F^*_\ell(s,\alpha)}{Q_\ell(s)} = \Gamma(2s-1/2) H(s),
\end{equation}
say. Moreover, thanks to the definition of $F_\ell^\pm(s,\alpha)$ and of $\nu_\pm$, we have
\begin{equation}
\label{3-2}
F_\ell^\pm(s,\alpha) = \frac{c^*(\nu_\pm^2)}{|\nu_\pm^2|^{1/2} |\nu_\pm+\nu_\alpha|^{2s-1/2}} \left(\frac{|\nu_\pm+\nu_\alpha|}{|\nu_\pm|}\right)^\ell (1+h_\ell^\pm(s)),
\end{equation}
where $h_\ell^\pm(s)$ are holomorphic and satisfy, uniformly for $t\in\RR$ and $\ell=0,\dots,h^*$ as $\si\to+\infty$,
\[
h_\ell^\pm(s) =o(1).
\]
Hence, recalling the definition of $m_\pm$ and $c^*_\pm$ in \eqref{1-19}, and that $\deg Q_\ell\geq 1$ for $\ell\geq 1$, from \eqref{1-12}, \eqref{3-1} and \eqref{3-2} we obtain
\begin{equation}
\label{3-3}
\Sigma(s) = \Gamma(2s-1/2) \left(e^{-i\pi s} \frac{c^*_+}{m_+^{2s}} (1+h^+(s)) + e^{i\pi s} \frac{c^*_-}{m_-^{2s}} (1+h^-(s)) \right),
\end{equation}
where $h^\pm(s)$ are holomorphic and satisfy, uniformly for $t\in\RR$ as $\si\to+\infty$,
\begin{equation}
\label{3-300}
h^\pm(s) =o(1).
\end{equation}

\smallskip
In order to detect the zeros of $\Sigma(1-s)$ with $\si<0$ we write
\[
W(s) = e^{-i\pi s} \frac{c^*_+}{m_+^{2s}} + e^{i\pi s} \frac{c^*_-}{m_-^{2s}} 
\]
and study the zeros of $W(1-s)$; an application of Rouch\'e's theorem will then allow to get corresponding results for $\Sigma(1-s)$. Since $-1=e^{i\pi}$, in view of \eqref{1-19} we have that $W(1-s)=0$ if and only if
\[
e^{-i\pi(1-\si-it)}\rho_+e^{i\theta_+} m_+^{-2(1-\si-it)} = e^{i\pi(1-\si-it)+i\pi}\rho_-e^{i\theta_-} m_-^{-2(1-\si-it)},
\]
i.e.
\begin{equation}
\label{3-4}
e^{i(\pi(\si-1)+\theta_++2t\log m_+)} e^{-\pi t+\log\rho_++2(\si-1)\log m_+} = e^{i(-\pi\si+\theta_-+2t\log m_-)} e^{\pi t+\log\rho_-+2(\si-1)\log m_-}.
\end{equation}
Hence the moduli of the two sides of \eqref{3-4} are equal provided
\begin{equation}
\label{3-5}
\frac{\si}{\pi}\log\left(\frac{m_+}{m_-}\right) -t + \frac{1}{2\pi} \log\left(\frac{\rho_+m_-^2}{\rho_-m_+^2}\right) =0,
\end{equation}
while the arguments are equal provided for some $k\in\ZZ$
\begin{equation}
\label{3-6}
2\pi\si -2t\log\left(\frac{m_-}{m_+}\right) +\theta_+-\theta_- - (2k+1)\pi = 0.
\end{equation}
Since the lines in \eqref{3-5} and \eqref{3-6} are orthogonal, as $k$ varies over $\ZZ$ they have infinitely many intersections  in the half-plane $\si<0$. Hence $W(1-s)$ has infinitely many zeros on the part of the line \eqref{3-5} with $\si<0$, and the line \eqref{3-5} is exactly the same line in the definition of $\LL_\epsilon(\alpha)$, see \eqref{1-20}. Moreover, the number of such zeros in $-R\leq \si<-\si_0$ clearly equals
\begin{equation}
\label{3-60}
c_2(\alpha)R +O_{\si_0}(1),
\end{equation}
with a certain $c_2(\alpha)\neq0$, for any $\si_0>0$.

\smallskip
Now we write
\[
V(s) = e^{-i\pi s} \frac{c^*_+}{m_+^{2s}} h^+(s) + e^{i\pi s} \frac{c^*_-}{m_-^{2s}} h^-(s),
\]
thus by \eqref{3-1} and \eqref{3-3}
\begin{equation}
\label{3-8}
H(1-s) = W(1-s) + V(1-s).
\end{equation}
In view of the above argument, Corollary 2 will follow if we show that there exists $c>0$ such that for every sufficiently small $\epsilon>0$
\begin{equation}
\label{3-9}
|W(1-s)| \geq c \epsilon \max \left(\left|e^{-i\pi(1-s)} \frac{c^*_+}{m_+^{2(1-s)}} \right|, \left|e^{i\pi(1-s)} \frac{c^*_-}{m_-^{2(1-s)}}\right|\right) = c\epsilon \rho(1-s),
\end{equation}
say, for $s\not\in\LL_\epsilon(\alpha)$ with $\si<0$ and for $s$ on the vertical segments $\R_n= \LL_\epsilon(\alpha)\cap \{\si=-R_n\}$, where $R_n\geq0$, $n\geq1$, is a suitable increasing sequence tending to $+\infty$. Indeed, by \eqref{3-300}, for $\si\to-\infty$ we have
\begin{equation}
\label{3-10}
V(1-s) = o\big(\rho(1-s)\big)
\end{equation}
uniformly for $t\in\RR$, hence \eqref{3-8}, \eqref{3-9} and \eqref{3-10} imply that $H(1-s)\neq0$ for $s\not\in\LL_\epsilon(\alpha)$ and $\si<-\si_\epsilon$ for some $\si_\epsilon\geq0$. Therefore, {\sl i)} of Corollary 2 follows from \eqref{1-16} and \eqref{3-1}. Moreover, given $n_0$ with $R_{n_0}>\si_\epsilon$, let $\P_\epsilon$ be the parallelogram with sides $\R_{n_0}$, $\R_n$ with $n>n_0$ and the two segments of boundary of $\LL_\epsilon(\alpha)$ with real part between $-R_n$ and $-R_{n_0}$. Then by the same reason as before, for $s\in\P_\epsilon$ we have
\[
|W(1-s)|>|V(1-s)|,
\]
hence by Rouch\'e's theorem $H(1-s)$ has inside $\P_\epsilon$ the same number of zeros as $W(1-s)$. Therefore, {\sl ii)} of Corollary 2 follows by letting $n\to\infty$, since we already detected the zeros of $W(1-s)$, and the poles of $\Gamma(2(1-s)-1/2)$ cannot cancel such zeros.

\smallskip
It remains to prove \eqref{3-9}. We already know that
\[
\left|e^{-i\pi(1-s)} \frac{c^*_+}{m_+^{2(1-s)}} \right|= \left|e^{i\pi(1-s)} \frac{c^*_-}{m_-^{2(1-s)}}\right| = \widetilde{\rho}(\si),
\]
say, when $s$ is on the line \eqref{3-5}. Writing such a line as $t=t(\si)$, we denote the points $s\not\in\LL_\epsilon(\alpha)$ with $\si<0$ as $s= \si+i(t(\si) +\delta)$ with $|\delta|\geq \epsilon$. Then in view of \eqref{3-4} we have
\[
W(1-s) = e^{i\theta_1(s)} \widetilde{\rho}(\si) e^{-\pi\delta} + e^{i\theta_2(s)} \widetilde{\rho}(\si) e^{\pi\delta}
\]
with certain $\theta_j(s)\in\RR$. Hence with a certain $c>0$ we have
\[
|W(1-s)| \geq \widetilde{\rho}(\si)\big(e^{\pi|\delta|} - e^{-\pi|\delta|}\big) \geq c\epsilon  \widetilde{\rho}(\si) e^{\pi|\delta|} =  c\epsilon \rho(1-s),
\]
as it is easy to check analysing two cases $\epsilon\leq |\delta|\leq 1$ and $|\delta|\geq 1$. A similar lower bound can be obtained for $s\in\R_n$ as follows. By continuity, between two consecutive intersections of the lines \eqref{3-5} and \eqref{3-6} in the half-plane $\si<0$, there is a value of $\si$ such that the two terms in $W(1-s)$ have the same argument when $s=\si+it(\si)$. Then we choose the $R_n$ to be these values of $\si$, thus
\[
\R_n=\{R_n+i(t(R_n) +\delta), |\delta|<\epsilon\}.
\]
Due to the shape of such arguments, see \eqref{3-4}, and since $0<\epsilon<c_1(\alpha)$ and $c_1(\alpha)$ can be chosen sufficiently small depending on $F(s)$ and $\alpha$, the absolute value of the difference between the above arguments when $s$ runs over $\R_n$ is bounded by $\pi/100$, say. Hence
\[
|W(1-s)| \geq c \rho(1-s)
\]
for $s\in\R_n$, and \eqref{3-9} follows. Corollary 2 is therefore proved. \fine

\bigskip
{\bf 3.2.~Proof of Corollary 3.} Let $N_F(T,\alpha)$ be as in \eqref{1-22}, $T_0>0, a>\si^-$ and $b>\si^+$ be sufficiently large,
\[
N^+(T) = \{\rho=\beta+i\gamma: F(\rho,\alpha)=0 \ \text{with} \ -a\leq \si\leq b \ \text{and} \ T_0<\gamma\leq T\}
\]
and similarly for $N^-(T)$, with $-T\leq \gamma<-T_0$. Then in view of the results in Subsection 1.6, see in particular \eqref{1-21}, we have
\begin{equation}
\label{3-12}
N_F(T,\alpha) = N^+(T) + N^-(T) + O(1).
\end{equation}
Corollary 3 will follow from a suitable application to $N^\pm(T)$ of the classical technique for the Riemann-von Mangoldt formula based on the argument principle; we give only a sketch of proof, mainly in order to compute the slightly nonstandard coefficient of $T$ in \eqref{1-23}.

\smallskip
Let $\R^+$ be the rectangle joining the points $-a+iT_0$, $b+iT_0$, $b+iT$ and $-a+iT$ with positive orientation, and denote by $\LL_j$, $j=1,\dots,4$, its sides, starting with the lower horizontal side. Hence
\begin{equation}
\label{3-13}
N^+(T) = \frac{1}{2\pi} \sum_{j=1}^4\Delta_{\LL_j} \arg F(s,\alpha) = \frac{1}{2\pi} (\Delta_1+\cdots+\Delta_4),
\end{equation}
say, and clearly
\begin{equation}
\label{3-14}
\Delta_1 = O(1).
\end{equation}
Recalling the definition of $\overline{n}$ in Subsection 1.7 we have
\[
F(s,\alpha) = \frac{1}{\overline{n}^s} \left(a(\overline{n})e(-\alpha\sqrt{\overline{n}}) + \sum_{n>\overline{n}} \frac{a(n)e(-\alpha\sqrt{n})}{(n/\overline{n})^s}\right) = \frac{1}{\overline{n}^s} \big(a(\overline{n})e(-\alpha\sqrt{\overline{n}})+ o(1)\big)
\]
uniformly for $s\in\LL_2$ as $b\to\infty$, hence for $b$ sufficiently large
\begin{equation}
\label{3-15}
\Delta_2 = -T\log \overline{n} + O(1).
\end{equation}
For the treatment of $\LL_3$ we consider the conjugate function
\[
\overline{F}(s,-\alpha) = \sum_{n=1}^\infty \frac{\overline{a(n)}e(\alpha\sqrt{n})}{n^s},
\]
so that
\[
2\Re{F(\si+iT,\alpha)} = F(\si+iT,\alpha) + \overline{F}(\si-iT,-\alpha) = f(\si),
\]
say. We follow the standard approach of bounding the variation of the argument of $F(s,\alpha)$ on $\LL_3$ by means of the number of zeros of $\Re{F(\si+iT,\alpha)}$ on the segment $-a\leq \si\leq b$. In turn, such number is bounded via Jensen's inequality by the number of zeros of the holomorphic function $f(s)$ on the circle of radius $r=(a+b)/2$ and center $s_0=(b-a)/2$. Since $F(s,\alpha)$ has polynomial growth on vertical strips we obtain
\begin{equation}
\label{3-16}
\Delta_3 = O(\log T).
\end{equation}
On $\LL_4$ we use functional equation \eqref{1-16}. Using once again the factorial formula for the $\Gamma$ function as in \eqref{2-44} and \eqref{3-1}, in view of the definition of $F^*_\ell(s,\alpha)$ in \eqref{1-12} and of $m_\pm$ in \eqref{1-19} we get, uniformly for $s\in\LL_4$ as $a\to\infty$, that
\[
F(s,\alpha) = c e^{-i\pi s}\left(\frac{\sqrt{N}}{4\pi m_-}\right)^{-2s} \Gamma(2(1-s)-1/2) \big(1+o(1)\big)
\]
with a certain constant $c\in\CC$, $c\neq0$. Hence for $a$ sufficiently large, from Stirling's formula we obtain
\begin{equation}
\label{3-17}
\Delta_4 = 2T\log T + T\log\left(\frac{N}{(2\pi em_-)^2}\right) + O(1),
\end{equation}
therefore from \eqref{3-13}-\eqref{3-17} we finally deduce that
\begin{equation}
\label{3-18}
N^+(T) = \frac{1}{\pi}T\log T + \frac{T}{2\pi}\log\left(\frac{N}{\overline{n}(2\pi em_-)^2}\right) + O(\log T).
\end{equation}

\smallskip
A completely analogous argument, applied to the rectangle $\R^-$ joining the points $-a-iT$, $b-iT$, $b-iT_0$ and $-a-iT_0$, with positive orientation, shows that
\begin{equation}
\label{3-19}
N^-(T) = \frac{1}{\pi}T\log T + \frac{T}{2\pi}\log\left(\frac{N}{\overline{n}(2\pi em_+)^2}\right) + O(\log T).
\end{equation}
Corollary 3 follows now from \eqref{3-12}, \eqref{3-18} and \eqref{3-19}. \fine

\bigskip
\section{Appendix}

\smallskip
Here we prove an assertion made in Subsection 1.3, namely that if a function $F\in\S^\sharp$ has degree $\geq2$ and satisfies the Ramanujan conjecture, then its standard twist $F(s,\alpha)$ does not satisfy a functional equation of type \eqref{1-2}. Since the argument is similar to the proof of Theorem 1 in \cite{Ka-Pe/2014a}, which asserts the slightly weaker statement that $F(s,\alpha)$ does not belong to $\S^\sharp$ under the same hypotheses on $F(s)$, we only give a sketch of the proof.

\smallskip
{\bf Fact 1.} {\sl The standard twist $L(s,\beta)$ of $L(s)$ is meromorphic on $\CC$ with all poles (if any) on a certain horizontal left half-line, and away from its poles is $\ll e^{|s|^c}$ for some $c>0$ as $|s|\to\infty$. Moreover, if $\beta\in$ Spec$(L)$ then $L(s,\beta)$ has a pole at the right end of such half-line. Further, all nonlinear twists of $L(s)$ of type
\[
\sum_{n=1}^\infty \frac{a(n)}{n^s} e(-\alpha_0 n^{1/d} - \alpha_1 n^\lambda),
\]
where $\alpha_0\geq 0$, $\alpha_1>0$ and $0<\lambda<1/d$, are entire.}

\smallskip
Fact 1 follows as a special case from Theorems 1, 2 and 3 of \cite{Ka-Pe/resoI}.

\smallskip
{\bf Fact 2.} {\sl If the standard twist $L(s,\alpha)$ satisfies a functional equation of type \eqref{1-2}, then the results in Fact $1$ hold also with $L(s,\alpha)$ in place of $L(s)$.}

\smallskip
Indeed, Theorems 1, 2 and 3 of \cite{Ka-Pe/resoI} are proved for functions belonging to $\S^\sharp$, but under the assumption in Fact 2 we have:

\smallskip
(i) $L(s,\alpha)$ is absolutely convergent for $\si>1$;

(ii) $L(s,\alpha)$ has the meromorphic structure described in Fact 1;

(iii) $L(s,\alpha)$ satisfies a functional equation of type \eqref{1-2}.

\smallskip
\noindent
Fact 2 now follows observing that (i)-(iii) differ from the definition of $\S^\sharp$ only by condition (ii); moreover, one can check that the arguments leading to the above quoted theorems do not depend on the presence of a pole at $s=1$ only, as in the case of $\S^\sharp$, and still work if the function under consideration satisfies (ii). From Facts 1 and 2, and the arguments in Theorem 1 of \cite{Ka-Pe/2014a}, we now derive the following result, thus justifying the assertion made in Subsection 1.3.

\medskip
{\bf Proposition 4.1.}  {\sl If $F\in\S^\sharp$ has degree $\geq 2$ and satisfies the Ramanujan conjecture, then the standard twist $F(s,\alpha)$ does not satisfy a functional equation of type \eqref{1-2}.}

\medskip
Note that Proposition 4.1 in the case $\alpha\not\in$ Spec$(F)$ is an immediate consequence of the above mentioned Theorem 1 in \cite{Ka-Pe/2014a}. Indeed, in this case $F(s,\alpha)$ satisfies the first two conditions in the definition of $\S^\sharp$, hence the fact that $F(s,\alpha)$ does not belong to $\S^\sharp$ means simply that $F(s,\alpha)$ does not satisfy functional equation \eqref{1-2}.

\medskip
{\it Proof of Proposition 4.1.} Let $d_F$ and $a(n)$ denote, respectively, degree and Dirichlet coefficients of $F(s)$, write
\[
G(s) = F(s,\alpha)
\]
and assume that $G(s)$ satisfies \eqref{1-2}. Then we can define degree $d_G$ and spectrum Spec$(G)$ of $G(s)$ and, thanks to Fact 2, $G(s)$ has the properties described in Fact 1 for $L(s)$.

\smallskip
We may assume that $d_G\geq d_F$, otherwise we replace $F(s)$ by $\overline{G}(s)$ in what follows. If $d_G>d_F$ we choose $\beta\in$ Spec$(G)$ and note from Fact 1, applied to $G(s)$, that the standard twist $G(s,\beta)$ is not entire. Moreover, it can be written as
\[
G(s,\beta) = \sum_{n=1}^\infty \frac{a_F(n)}{n^s} e(-\alpha n^{1/d_F} - \beta n^{1/d_G}).
\]
But, since $\beta>0$ and $0<1/\d_G<1/d_F$, from Fact 1 applied to $F(s)$ we have that $G(s,\beta)$ is entire, a contradiction. 

\smallskip
Thus we have that $d_F=d_G=d\geq 2$, say. Now we follow the proof of Theorem 1 of \cite{Ka-Pe/2014a} (see p.150, after (2.1) there) and obtain that the Dirichlet series of $G(s)$, and hence of $F(s)$ as well, is absolutely convergent for $\si>1/d$. This gives a contradiction, since Corollary 3 of \cite{Ka-Pe/2005} states that the abscissa of absolute convergence of $F(s)$ is $\geq (d+1)/2d$. \fine

\medskip
We finally remark that the Ramanujan conjecture enters the part of the proof of Theorem 1 in \cite{Ka-Pe/2014a} which we followed without giving details, and is used to show the absolute convergence for $\si>1/d$.

\newpage

\ifx\undefined\bysame{poly}.
\newcommand{\bysame}{\leavevmode\hbox to3em{\hrulefill}\ ,}
\fi

\bigskip
\bigskip
\bigskip
\noindent
Jerzy Kaczorowski, Faculty of Mathematics and Computer Science, A.Mickiewicz University, 61-614 Pozna\'n, Poland and Institute of Mathematics of the Polish Academy of Sciences, 
00-956 Warsaw, Poland. e-mail: kjerzy@amu.edu.pl

\medskip
\noindent
Alberto Perelli, Dipartimento di Matematica, Universit\`a di Genova, via Dodecaneso 35, 16146 Genova, Italy. e-mail: perelli@dima.unige.it


\begin{thebibliography}{100} {\normalsize

\bibitem{Ber/1972} B.C.Berndt - {\sl Two new proofs of Lerch's functional equation} - Proc. A. M. S. {\bf 32} (1972), 403--408.

\bibitem{Bru/1997} J.H.Bruinier - {\sl Modulformen halbganzen Gewichts und Beziehung zu Dirichletreihen} - Diplomarbeit, Universit\"at Heidelberg, 1997 (\url{http://www.mathematik.tu-darmstadt.de/fbereiche/AlgGeoFA/staff/bruinier/publications/dipl.pdf}).

\bibitem{C-M-P/2009} E.Carletti, G.Monti Bragadin, A.Perelli - {\sl A note on Hecke's functional equation and the Selberg class} - Funct. Approx. {\bf 41} (2009), 211--220.

\bibitem{Co-Gh/1993} J.B.Conrey, A.Ghosh - {\sl On the Selberg class of Dirichlet series: small degrees} - Duke Math. J. {\bf 72} (1993),  673--693.

\bibitem{Iwa/1997} H.Iwaniec - {\sl Topics in Classical Automorphic Forms} - A.M.S. Publications 1997.

\bibitem{Kac/2004} J.Kaczorowski -  {\sl Some remarks on Fourier coefficients of Hecke modular functions} - Comment.
Math., special volume {\it in honorem} J. Musielk, (2004), 105--121.

\bibitem{Kac/2006} J.Kaczorowski - {\sl Axiomatic theory of $L$-functions: the Selberg class} - In {\sl Analytic Number Theory}, C.I.M.E. Summer School, Cetraro (Italy), ed. by A.Perelli and C.Viola, 133--209, Springer L.N. 1891, 2006.

\bibitem{KMPSW/2006} J.Kaczorowski, G.Molteni, A.Perelli, J.Steuding, J.Wolfart - {\sl Hecke's theory and the Selberg class} - Funct. Approx. {\bf 35} (2006), 183--193.

\bibitem{Ka-Pe/1999a} J.Kaczorowski, A.Perelli - {\sl On the structure of the Selberg class, I: $0\leq d \leq 1$} - Acta Math. {\bf 182} (1999), 207--241.

\bibitem{Ka-Pe/1999b} J.Kaczorowski, A.Perelli - {\sl The Selberg class: a survey} - In {\sl Number Theory in  Progress}, Proc. Conf. in Honor of A.Schinzel, ed. by K.Gy\"ory {\sl et al.},  953--992, de Gruyter 1999.

\bibitem{Ka-Pe/2005} J.Kaczorowski, A.Perelli - {\sl On the structure of the Selberg class, VI: non-linear twists} - Acta Arith. {\bf 116} (2005), 315--341.

\bibitem{Ka-Pe/2011} J.Kaczorowski, A.Perelli - {\sl On the structure of the Selberg class, VII: $1<d<2$} - Annals of Math. {\bf 173} (2011), 1397--1441.

\bibitem{Ka-Pe/2014a} J.Kaczorowski, A.Perelli -  {\sl Internal twists of $L$-functions} - In {\sl Number Theory, Analysis, and Combinatorics}, ed. by J.Pintz {\it et al.}, 145--154, de Gruyter 2014.

\bibitem{Ka-Pe/resoI} J.Kaczorowski, A.Perelli - {\sl Twists and resonance of $L$-functions, I} - J. European Math. Soc. {\bf 18} (2016), 1349--1389.

\bibitem{Ka-Pe/abs} J.Kaczorowski, A.Perelli - {\sl Some remarks on the convergence of the Dirichlet series of $L$-functions and related questions} - Math. Zeitschrift {\bf 285} (2017), 1345--1355.

\bibitem{Ka-Pe/2017} J.Kaczorowski, A.Perelli - {\sl A note on Linnik's approach to the Dirichlet $L$-functions} - Trudy Mat. Inst. Steklova  {\bf 296} (2017), 123--132 (Russian); English transl. Proc. Steklov Inst. Math. {\bf 296} (2017), 115-124.

\bibitem{Ka-Pe/book} J.Kaczorowski, A.Perelli - {\sl Introduction to the Selberg Class of $L$-Functions} - In preparation.

\bibitem{Miy/1989} T.Miyake - {\sl Modular Forms} - Springer Verlag 1989.

\bibitem{Ogg/1969} A.Ogg - {\sl Modular Forms and Dirichlet Series} - Benjamin 1969.

\bibitem{Pa-Ka/2001} R.B.Paris, D.Kaminski - {\sl Asymptotics and Mellin-Barnes Integrals} - Cambridge U. P. 2001.

\bibitem{Per/2005} A.Perelli - {\sl A survey of the Selberg class of $L$-functions, part I} - Milan J. Math. {\bf 73} (2005), 19--52.

\bibitem{Per/2004} A.Perelli - {\sl A survey of the Selberg class of $L$-functions, part II} - Riv. Mat. Univ. Parma (7) {\bf 3*} (2004), 83--118.

\bibitem{Per/2010} A.Perelli - {\sl Non-linear twists of $L$-functions: a survey} - Milan J. Math. {\bf 78} (2010), 117--134.

\bibitem{Per/2017} A.Perelli - {\sl Converse theorems: from the Riemann zeta function to the Selberg class} - Bollettino U.M.I. {\bf 10} (2017), 29--53.

\bibitem{Sel/1989} A.Selberg - {\sl Old and new conjectures and results about a class of Dirichlet series} - In {\sl Proc. Amalfi Conf. Analytic Number Theory}, ed. by E.Bombieri {\sl et al.}, 367--385, Universit\`a di Salerno 1992; {\sl Collected Papers}, vol. II, 47--63, Springer Verlag 1991.

\bibitem{Tit/1939} E.C.Titchmarsh - {\sl The Theory of Functions} - $2^{\text{nd}}$ ed., Oxford U. P. 1939.

}
\end{thebibliography}
\end{document}